\documentclass[12pt]{article}
\usepackage{latexsym} \usepackage{epsf}
\usepackage{a4}
\usepackage{amsfonts}
\renewcommand{\theequation}{\thesection.\arabic{equation}}
\newcounter{subequation}[equation]
\makeatletter

\expandafter\let\expandafter\reset@font\csname reset@font\endcsname

\def\subeqnarray{\arraycolsep1pt
    \def\@eqnnum\stepcounter##1{\stepcounter{subequation}%
        {\reset@font\rm(\theequation\alph{subequation})}}
\jot5mm     \eqnarray}

\makeatother

\def\be{\begin{equation}}
\def\ee{\end{equation}}
\def\bea{\begin{eqnarray}}
\def\eea{\end{eqnarray}}
\def\dd{\partial}
\def\half{\frac{1}{2}}

\def\one#1{#1^{\raise5pt\hbox{$\scriptstyle\!\!\!\!1$}}\,{}}
\def\two#1{#1^{\raise5pt\hbox{$\scriptstyle\!\!\!\!2$}}\,{}}

\def\II{\hbox{{1}\kern-.25em\hbox{l}}}
\makeatletter
\def\binrel@#1{\begingroup
  \setboxz@h{\thinmuskip0mu
    \medmuskip\m@ne mu\thickmuskip\@ne mu
    \setbox\tw@\hbox{$#1\m@th$}\kern-\wd\tw@
    ${}#1{}\m@th$}%
  \edef\@tempa{\endgroup\let\noexpand\binrel@@
    \ifdim\wdz@<\z@ \mathbin
    \else\ifdim\wdz@>\z@ \mathrel
    \else \relax\fi\fi}%
  \@tempa
}
\let\binrel@@\relax
\def\overset#1#2{\binrel@{#2}%
  \binrel@@{\mathop{\kern\z@#2}\limits^{#1}}}
\def\underset#1#2{\binrel@{#2}%
  \binrel@@{\mathop{\kern\z@#2}\limits_{#1}}}
\makeatother
\newfont{\bbd}{msbm10 scaled\magstep1}
\def\C{\hbox{\bbd C}}

\def\F{{\mathcal R}}

\def\R{\hbox{\bbd R}}

\def\S{\hbox{\bbd S}}

\def\P{\hbox{\bbd P}}

\def\RR{{\mathcal R}}

\newtheorem{prop}{Proposition}
\newtheorem{lem}{Lemma}

\begin{document}

\vskip 1cm

\centerline{\LARGE\bf Factorization of the R-matrix. II.}

\vskip 1cm \centerline{\sc S.E. Derkachov} \vskip 1cm

\centerline{Department of Mathematics, St Petersburg
Technology Institute}\centerline{ St.Petersburg, Russia.}
\centerline{E-mail: {\tt S.Derkachov@pobox.spbu.ru}}

\vskip 2cm

{\bf Abstract.} We study the general rational solution of
the Yang-Baxter equation with the supersymmetry algebra
$s\ell(2|1)$. The R-operator acting in the tensor product
of two arbitrary representations of the supersymmetry
algebra can be represented as the product of the simpler
"building blocks" -- $\RR$-operators.

\newpage

{\small \tableofcontents}
\renewcommand{\refname}{References.}
\renewcommand{\thefootnote}{\arabic{footnote}}
\setcounter{footnote}{0}
\setcounter{equation}{0}

\section{Introduction}
\setcounter{equation}{0}

In the previous paper~\cite{I} we have shown that the
general R-matrix can be represented as the product of the
much more simpler $\F$-operators. In this paper we shall
consider the general rational solution of the Yang-Baxter
equation with the supersymmetry algebra $s\ell(2|1)$ and
show that there exists the same factorization. In fact all
calculations are very similar to the $s\ell(3)$-example and
modifications due to supersymmetry are simple. The
generalization of the previous results~\cite{I} to the
algebra of supersymmetry is mainly motivated by the
possible applications to the super Yang-Mills
theory~\cite{B,DNW,K,BDKM}.

The presentation is organized as follows. In Section 2 we
collect the standard facts about the algebra~$s\ell(2|1)$
and its representations. We represent the lowest weight
modules by polynomials in one even variable~($z$) and two
odd variables~($\theta,\bar\theta$) and the
$s\ell(2|1)$-generators as first order differential
operators. We use the notation and formulae from the
paper~\cite{DKK}. In Section 3 we derive the defining
relation for the general R-matrix, i.e. the solution of the
Yang-Baxter equation acting on tensor products of two
arbitrary representations, the elements of which are
polynomials in variables~$z_1,\theta_1,\bar\theta_1$
and~$z_2,\theta_2,\bar\theta_2$. In Section 4 we introduce
the natural defining equations for the $\F$-operators and
show that the general R-matrix can be represented as the
product of such much more simple operators.

Finally, in Section 4 we summarize. In Appendix we
calculate the matrix elements of the $\F$-operators and as
consequence obtain the matrix elements of R-matrix in full
agreement with the results of the paper~\cite{DKK}.

\section{$sl(2|1)$ lowest weight modules}
\setcounter{equation}{0}

The superalgebra $s\ell(2|1)$ has eight generators: four
odd $\mathbf{V}^{\pm},\mathbf{W}^{\pm}$ and four even
$\mathbf{S},\mathbf{S}^{\pm}$ and $\mathbf{B}$. Using the
natural notations
$$
\mathbf{E}_{31} = \mathbf{S}^{-}\ ;\ \mathbf{E}_{21} = -
\mathbf{W}^{-}\ ;\ \mathbf{E}_{32} = \mathbf{V}^{-}\ ;\
\mathbf{E}_{13} = \mathbf{S}^{+}\ ;\ \mathbf{E}_{23} =
\mathbf{W}^{+}\ ;\ \mathbf{E}_{12} = \mathbf{V}^{+}
$$
$$
\mathbf{E}_{11} = \mathbf{B} - \mathbf{S}\ ;\
\mathbf{E}_{22} = -2 \mathbf{B}\ ;\ \mathbf{E}_{33} =
\mathbf{B} + \mathbf{S}.
$$
the commutation relations for the generators of
$s\ell(2|1)$ can be written compactly in the
form~\cite{SNR,JG}
$$
[\mathbf{E}_{AB},\mathbf{E}_{CD}] = \delta_{CB}
\mathbf{E}_{AD} - (-)^{(\bar{A}+\bar{B})(\bar{C}+\bar{D})}
\delta_{AD} \mathbf{E}_{CB}\  ;\ A,B,C,D = 1,2,3
$$
where the graded commutator is defined as~(we choose the
grading $\bar 1=\bar 3=0\ ,\ \bar 2 = 1$)
$$
[\mathbf{E}_{AB},\mathbf{E}_{CD}]\equiv
\mathbf{E}_{AB}\cdot \mathbf{E}_{CD} -
(-)^{(\bar{A}+\bar{B})(\bar{C}+\bar{D})}
\mathbf{E}_{CD}\cdot \mathbf{E}_{AB}.
$$
There are two central elements~\cite{SNR,JG,DIC}
$$
\mathbf{C}_2 = \frac{1}{2}\sum_{AB} (-)^{\bar B}
\mathbf{E}_{AB}\mathbf{E}_{BA} =
\mathbf{S}^2-\mathbf{B}^2+\mathbf{S}^{+}\mathbf{S}^{-}+
\mathbf{V}^{+}\mathbf{W}^{-}+\mathbf{W}^{+}\mathbf{V}^{-}\
;\ \mathbf{C}_3 = \frac{1}{6}\sum_{ABC} (-)^{\bar B +\bar
C} \mathbf{E}_{AB}\mathbf{E}_{BC}\mathbf{E}_{CA}
$$
The Verma module is the generic lowest weight
$s\ell(2|1)$-module $\mathbf{V}_{\Lambda}\ ;\ \Lambda =
(\ell , b) $. As a linear space $\mathbf{V}_{\Lambda}$ is
spanned by the basis with even elements $\mathbf{a}_k ,
\mathbf{b}_k$
$$
\mathbf{a}_k = \mathbf{S}_{+}^k \mathbf{a}_0 \ ,\ k=0 , 1 ,
2\cdots \ ;\ \mathbf{b}_k =
\mathbf{S}_{+}^{k-1}\mathbf{W}_{+}\mathbf{V}_{+}\mathbf{a}_0\
\ ,\ k= 1 , 2 \cdots
$$
and odd elements $\mathbf{v}_k , \mathbf{w}_k$
$$
\mathbf{v}_k = \mathbf{S}_{+}^k \mathbf{V}_{+}\mathbf{a}_0
\ ;\ \mathbf{w}_k =
\mathbf{S}_{+}^{k}\mathbf{W}_{+}\mathbf{a}_0\ ,\ k=0 , 1 ,
2\cdots
$$
The vector $\mathbf{a}_0$ is the lowest weight vector:
$$\mathbf{S}_{-}\mathbf{a}_0
=\mathbf{V}_{-}\mathbf{a}_0 = \mathbf{W}_{-}\mathbf{a}_0
 = 0 \ ;\ \mathbf{S}\mathbf{a}_0 = \ell\cdot\mathbf{a}_0
 \ ;\  \mathbf{B}\mathbf{a}_0 = b\cdot\mathbf{a}_0
$$
We shall use the representation $\mathrm{V}_{\Lambda}$ of
$s\ell(2|1)$ in the infinite-dimensional space $\C[Z]$
where $Z=(z,\theta,\bar\theta)$ of polynomials in even
variable $z$ and odd variables $\theta,\bar\theta$ with the
monomial basis $\left\{z^k , \theta\bar\theta z^k ;\ \theta
z^k , \bar\theta z^k\right\}_{k= 0}^{\infty}$ and lowest
weight vector $a_0 = 1$~\cite{DKK}. The action of
$s\ell(2|1)$ in $\mathrm{V}_{\Lambda}$ is given by the
first-order differential operators \be \mathrm{S}^{-} =
-\partial \ ;\ \mathrm{V}^{-} =
\partial_{\theta}+\frac{1}{2}\bar\theta\partial \ ;\ \mathrm{W}^{-}
= \partial_{\bar\theta}+\frac{1}{2}\theta\partial
\label{gen} \ee
$$
\mathrm{V}^{+}= -\left[z\partial_{\theta}
+\frac{1}{2}\bar\theta z\partial+
\frac{1}{2}\bar\theta\theta\partial_{\theta}\right]-
(\ell-b)\bar\theta \ ;\ \mathrm{W}^{+}=
-\left[z\partial_{\bar\theta} +\frac{1}{2}\theta z\partial+
\frac{1}{2}\theta\bar\theta\partial_{\bar\theta}\right]-
(\ell+b)\theta
$$
$$
\mathrm{S}^{+}= z^2\partial + z\theta\partial_{\theta}+
z\bar\theta\partial_{\bar\theta} + 2\ell z
-b\theta\bar\theta \ ;\ S= z\partial
+\frac{1}{2}\theta\partial_{\theta}+
\frac{1}{2}\bar\theta\partial_{\bar\theta} +\ell \ ;\  B=
\frac{1}{2}\bar\theta\partial_{\bar\theta} -
\frac{1}{2}\theta\partial_{\theta}+b
$$
It is possible to derive the closed expressions for the
elements of the basis
$$
a_{k}= \mathrm{S}^k_{+}\cdot 1 =(2\ell)_k \left[z- \frac{k
b}{2\ell}\cdot\theta\bar\theta\right]\cdot z^{k-1}\ ;\
b_{k} =
\mathrm{S}^{k-1}_{+}\mathrm{W}_{+}\mathrm{V}_{+}\cdot 1 =
\frac{\ell-b}{2\ell}(2\ell)_k\left[z+
\left(b+\ell+\frac{k}{2}\right)\cdot\theta\bar\theta\right]\cdot
z^{k-1}
$$
$$
v_{k}= \mathrm{S}^{k}_{+}\mathrm{V}_{+}\cdot 1 =
-(\ell-b)(2\ell+1)_k z^k\bar\theta \ ;\
w_{k}=\mathrm{S}^{k}_{+}\mathrm{W}_{+}\cdot 1 =
-(\ell+b)(2\ell+1)_k z^k\theta \ ;\ (2\ell)_k \equiv
\frac{\Gamma(2\ell+k)}{\Gamma(2\ell)}
$$
It is evident that for the generic $\ell\neq -\frac{n}{2}$
the module $\mathrm{V}_{\Lambda}$ is an irreducible lowest
weight $s\ell(2|1)$-module isomorphic to
$\mathbf{V}_{\Lambda}$ but for the special values of the
spin $\ell = - \frac{n}{2}$ there exists the finite
dimensional invariant subspace. There are three cases
depending on the relation between $b$ and
$n$~\cite{SNR,Mar,DIC}. The first case is for generic $b
\ne \pm \frac{n}{2}$(typical representations) and there
exists the $4n$-dimensional invariant subspace. The second
and third cases appear for~$b = \pm \frac{n}{2}$ (atypical
representations). For the chiral representation $b =
-\frac{n}{2}$ the $(2n+1)$-dimensional invariant subspace
is spanned on the vectors
$$
\Phi^{+}_k = \left(z - \frac{\theta\bar\theta}{2}\right)^k
\ ,\ k= 0...n \ ;\ \mathrm{W}_k = \theta z^{k}\ ,\ k= 0...n
- 1
$$
and for the antichiral representation $b = \frac{n}{2}$ the
$(2n+1)$-dimensional invariant subspace is spanned on the
vectors
$$
\Phi^{-}_k = \left(z +\frac{\theta\bar\theta}{2}\right)^k \
,\ k= 0...n \ ;\ \mathrm{V}_k = \bar\theta z^{k}\ ,\ k=
0...n - 1
$$
We shall use the three-dimensional chiral representation
$\mathrm{V}$.In the basis
$$
\mathbf{e}_1 = \mathrm{S}_{+}\cdot 1 =
-z+\frac{\theta\bar\theta}{2}\ ,\ \mathbf{e}_2 =
\mathrm{W}_{+}\cdot 1 = \theta\ ,\ \mathbf{e}_3 = 1
$$
the $s\ell(2|1)$-generators take the form \be
\mathbf{s}_{-} = \left (\begin{array}{ccc}
0 & 0 & 0  \\
0 & 0 & 0  \\
1 & 0 & 0
\end{array} \right )
\ ;\  \mathbf{w}_{-} = \left (\begin{array}{ccc}
0 & 0 & 0  \\
-1 & 0 & 0  \\
0 & 0 & 0
\end{array} \right )
\ ;\  \mathbf{v}_{-} = \left (\begin{array}{ccc}
0 & 0 & 0  \\
0 & 0 & 0  \\
0 & 1 & 0
\end{array} \right )\ ;\ \mathbf{s} = \left (\begin{array}{ccc}
\half & 0 & 0  \\
0 & 0 & 0  \\
0 & 0 & -\half
\end{array} \right )
\label{sfun} \ee
$$
\mathbf{s}_{+} = \left (\begin{array}{ccc}
0 & 0 & 1  \\
0 & 0 & 0  \\
0 & 0 & 0
\end{array} \right )
\ ;\  \mathbf{w}_{+} = \left (\begin{array}{ccc}
0 & 0 & 0  \\
0 & 0 & 1  \\
0 & 0 & 0
\end{array} \right )
\ ;\  \mathbf{v}_{+} = \left (\begin{array}{ccc}
0 & 1 & 0  \\
0 & 0 & 0  \\
0 & 0 & 0
\end{array} \right )
\ ;\  \mathbf{b} = \left (\begin{array}{ccc}
-\half & 0 & 0  \\
0 & -1 & 0  \\
0 & 0 & -\half
\end{array} \right )
$$
There exists the second three-dimensional representation -
antichiral representation $\bar{\mathrm{V}}$.  In the basis
$$
\mathbf{e}_1 = \mathrm{S}_{+}\cdot 1 =
-z-\frac{\theta\bar\theta}{2}\ ,\ \mathbf{e}_2 =
\mathrm{V}_{+}\cdot 1 = \bar\theta\ ,\ \mathbf{e}_3  = 1
$$
the $s\ell(2|1)$-generators take the similar form as in
chiral representation but $ \mathbf{v}^{\pm}\leftrightarrow
\mathbf{w}^{\pm}\ ;\ \mathbf{b}\rightarrow - \mathbf{b}$.
We use the standard definition for the matrix
$\mathrm{A}_{ik}$ of the linear operator $\mathbf{A}$ in
the basis $\{\mathbf{e}_k\}$
$$
\mathbf{A} \mathbf{e}_k = \sum_{i} \mathbf{e}_i
\mathrm{A}_{ik}
$$

\section{Yang-Baxter equation and Lax operator}
\setcounter{equation}{0}

The Yang-Baxter equation is the following three term
relation~\cite{KS,K}
$$
\R_{\Lambda_1\Lambda_2}(u-v)\R_{\Lambda_1\Lambda_3}(u)
\R_{\Lambda_2\Lambda_3}(v)=
\R_{\Lambda_2\Lambda_3}(v)\R_{\Lambda_1\Lambda_3}(u)
\R_{\Lambda_1\Lambda_2}(u-v)
$$
for the operators $\R_{\Lambda_i\Lambda_j}(u):
V_{\Lambda_i}\otimes V_{\Lambda_j}\to V_{\Lambda_i}\otimes
V_{\Lambda_j}$. We start from the simplest solutions of
Yang-Baxter equation and derive the defining equation for
the general $\R$-operator~\cite{DKK,KRS}. First we put
$\Lambda_1 = \Lambda_2 = \Lambda_3 =
(-\frac{1}{2},-\frac{1}{2})$ in Yang-Baxter equation and
consider the restriction on the invariant subspace
$\mathrm{V}\otimes\mathrm{V}\otimes\mathrm{V}$. We obtain
the equation
$$ \R_{1 2}(u-v)\R_{1 3}(u) \R_{2 3}(v)= \R_{2
3}(v)\R_{1 3}(u) \R_{1 2}(u-v) $$ where the operator
$\R_{12}(u)$ acts on the first and second copy of
$\mathrm{V}$ in the tensor product
$\mathrm{V}\otimes\mathrm{V}\otimes\mathrm{V}$ and
similarly for the other $\R$-operators. The solution is
well known~\cite{KS,K}
$$
\R_{12}(u) = u + \mathrm{P}_{12}
$$
where $\mathrm{P}_{12}$ is the (graded)permutation operator
in $\mathrm{V}\otimes\mathrm{V}$. We choose the basis
$\{\mathbf{e}_1,\mathbf{e}_2,\mathbf{e}_3\}$ in
$\mathrm{V}$ so that $\mathbf{e}_1,\mathbf{e}_3$ are even
elements and $\mathbf{e}_2$ is odd element and our grading
is: $\bar{1}=\bar{3}=0\ ,\ \bar{2}=1$. The permutation
operator acts on the basis as follows
$$
\mathrm{P}_{12}\mathbf{e}_i\otimes\mathbf{e}_k =
(-1)^{\bar{i}\bar{k}} \mathbf{e}_k\otimes\mathbf{e}_i
$$
and additional sign arises for
$\mathbf{e}_2\otimes\mathbf{e}_2$ only. Secondly we choose
$\Lambda_1 = \Lambda_2 = (-\frac{1}{2},-\frac{1}{2})\ ;\
\Lambda_3 = \Lambda = (\ell,b)$ and consider the
restriction on the invariant subspace
$\mathrm{V}\otimes\mathrm{V}\otimes\mathrm{V}_{\Lambda}$.
The restriction of the operator $\R_{\Lambda_1 \Lambda}(u)$
to the space $\mathrm{V}\otimes\mathrm{V}_{\Lambda}$
coincides up to normalization and shift of spectral
parameter with the Lax-operator~\cite{KS,K,rep}
$$
\mathrm{L}(u): \mathrm{V}\otimes\mathrm{V}_{\Lambda}\to
\mathrm{V}\otimes\mathrm{V}_{\Lambda}
$$
and the Yang-Baxter equation coincides with the defining
equation for the Lax-operator
$$
\R_{12}(u-v)\mathrm{L}^{(1)}(u)\mathrm{L}^{(2)}(v) =
\mathrm{L}^{(2)}(v)\mathrm{L}^{(1)}(u)\R_{12}(u-v)
$$
where $\mathrm{L}^{(1)}(u)$ is the operator which acts
nontrivially on the first copy of $\mathrm{V}$ and
$\mathrm{V}_{\Lambda}$ in the tensor product
$\mathrm{V}\otimes\mathrm{V}\otimes\mathrm{V}_{\Lambda}$
and $\mathrm{L}^{(2)}(u)$ is the operator which acts
nontrivially on the second copy of $\mathrm{V}$ and
$\mathrm{V}_{\Lambda}$. The solution coincides up to
additive constant with the Casimir operator $\mathbf{C}_2$
for the representation
$\mathrm{V}\otimes\mathrm{V}_{\Lambda}$~\cite{K,FPT}
$$
\mathrm{L}(u) \equiv  u  +
2\mathbf{s}\otimes\mathrm{S}-2\mathbf{b}\otimes\mathrm{B}+
\mathbf{v}_{+}\otimes
\mathrm{W}_{-}+\mathbf{s}_{+}\otimes\mathrm{S}_{-}
-\mathbf{w}_{-}\otimes\mathrm{V}_{+} +\mathbf{w}_{+}\otimes
\mathrm{V}_{-} +\mathbf{s}_{-}\otimes\mathrm{S}_{+}-
\mathbf{v}_{-}\otimes\mathrm{W}_{+}
$$
where $\mathbf{s} , \mathbf{b} , \mathbf{s}_{\pm} ,
\mathbf{v}_{\pm} , \mathbf{w}_{\pm}$ are
$s\ell(2|1)$-generators in the chiral
representation~(\ref{sfun}) and $\mathrm{S} , \mathrm{B} ,
\mathrm{S}_{\pm} , \mathrm{V}_{\pm} , \mathrm{W}_{\pm}$ are
generators in the generic representation. The algebra
$s\ell(2|1)$ has two three-dimensional representations --
chiral $\mathrm{V}$ and antichiral $\bar{\mathrm{V}}$ so
that there exists the second Lax-operator
$$
\bar{\mathrm{L}}(u):
\bar{\mathrm{V}}\otimes\mathrm{V}_{\Lambda}\to
\bar{\mathrm{V}}\otimes\mathrm{V}_{\Lambda}\ ;\
\R_{12}(u-v)\bar{\mathrm{L}}^{(1)}(u)\bar{\mathrm{L}}^{(2)}(v)
=
\bar{\mathrm{L}}^{(2)}(v)\bar{\mathrm{L}}^{(1)}(u)\R_{12}(u-v)
$$
The explicit expression for the second Lax-operator is the
same but now $\mathbf{s} , \mathbf{b} , \mathbf{s}_{\pm} ,
\mathbf{v}_{\pm} , \mathbf{w}_{\pm}$ are
$s\ell(2|1)$-generators in the antichiral representation
$\bar{\mathrm{V}}$. In matrix form we obtain~\cite{FPT,DKK}
$$
\mathrm{L}(u) = \left(\begin{array}{ccc}
\mathrm{S}+\mathrm{B}+u & -\mathrm{W}_{-} & \mathrm{S}_{-}\\
\mathrm{V}_{+} & 2\mathrm{B}+u & \mathrm{V}_{-}\\
\mathrm{S}_{+} & \mathrm{W}_{+} & \mathrm{B}-\mathrm{S}+u
   \end{array}\right)\ ;\ \bar{\mathrm{L}}(u) =
   \left(\begin{array}{ccc}
\mathrm{S}-\mathrm{B}+u & -\mathrm{V}_{-} & \mathrm{S}_{-}\\
\mathrm{W}_{+} & -2\mathrm{B}+u & \mathrm{W}_{-}\\
\mathrm{S}_{+} & \mathrm{V}_{+} & -\mathrm{B}-\mathrm{S}+u
   \end{array}\right)
$$
Note the arising some additional signs due to grading. For
example, we have $\mathbf{v}_{+}\mathbf{e}_2 =
\mathbf{e}_1$ but
$\mathbf{v}_{+}\otimes\mathrm{W}_{-}\mathbf{e}_2 =
-\mathbf{v}_{+}\mathbf{e}_2\mathrm{W}_{-} = -
\mathbf{e}_1\mathrm{W}_{-}$ so that one obtains
$$
\mathbf{v}_{+}\mathbf{e}_2 = \mathbf{e}_1\Rightarrow
\mathbf{v}_{+} = \left (\begin{array}{ccc}
0 & 1 & 0  \\
0 & 0 & 0  \\
0 & 0 & 0
\end{array} \right )\ ;\
\mathbf{v}_{+}\otimes\mathrm{W}_{-}\mathbf{e}_2 =
-\mathbf{e}_1\mathrm{W}_{-}\Rightarrow
\mathbf{v}_{+}\otimes\mathrm{W}_{-} = \left
(\begin{array}{ccc}
0 & -\mathrm{W}_{-}& 0  \\
0 & 0 & 0  \\
0 & 0 & 0
\end{array} \right )
$$
We shall use the chiral Lax-operator $\mathrm{L}(u)$ in
defining equation for the general $\R$-operator. The
Lax-operator depends on three parameters $u , \ell , b$ and
we shall use the parametrization
$$
u_1 = u+b+\ell\ ,\ u_2 = u+2b\ ,\ u_3 = u+b-\ell \ ;\ \ell
= \frac{u_1-u_3}{2}\ ,\ b = u_2-\frac{u_1+u_3}{2}
$$
The explicit form of the Lax operator
$\mathrm{L}(u_1,u_2,u_3)$ in the functional representation
$\mathrm{V}_{\Lambda}$ is \be \mathrm{L}(u_1,u_2,u_3) =
\left(\begin{array}{ccc} z\partial
+\bar\theta\partial_{\bar\theta}+u_1 &
-\left(\partial_{\bar\theta}+\frac{1}{2}\theta
\partial\right) &
-\partial\\
\mathrm{L}_{21}&
\bar\theta\partial_{\bar\theta}-\theta\partial_{\theta}+u_2
&
\partial_{\theta}+\frac{1}{2}\bar\theta\partial\\
\mathrm{L}_{31}& \mathrm{L}_{32} & -z\partial -
\theta\partial_{\theta}+u_3
   \end{array}\right)
\label{Lax21} \ee
$$
\mathrm{L}_{21} =
-\left(z-\frac{\theta\bar\theta}{2}\right)
\partial_{\theta}-\frac{1}{2}\bar\theta
z\partial + (u_2-u_1)\bar\theta\ ;\ \mathrm{L}_{32} =
-\left(z+\frac{\theta\bar\theta}{2}\right)
\partial_{\bar\theta}-\frac{1}{2}\theta z\partial
+(u_3-u_2)\theta
$$
$$
\mathrm{L}_{31} = z^2\partial
+z(\theta\partial_{\theta}+\bar\theta\partial_{\bar\theta})+
(u_1-u_3)z +(u_1+u_3-2u_2)\frac{\theta\bar\theta}{2}
$$
There exists the useful factorized representation for the
Lax-operator
\begin{equation}
\mathrm{L}(u_1,u_2,u_3) \equiv \left(\begin{array}{ccc}
1 & 0 & 0\\
-\bar{\theta} & 1 & 0\\
z+\frac{\theta\bar{\theta}}{2} & -\theta & 1
   \end{array}\right)
   \left(\begin{array}{ccc}
u_1 & \mathrm{D}^{-} & -\dd \\
0 & u_2-1 & -\mathrm{D}^{+}\\
0 & 0 & u_3
   \end{array}\right)
   \left(\begin{array}{ccc}
1 & 0 & 0\\
\bar{\theta} & 1 & 0\\
-z+\frac{\theta\bar{\theta}}{2} & \theta & 1
   \end{array}\right)
\label{factor21}
\end{equation}
where $\mathrm{D}^{\pm}$ are covariant derivatives
$$
\mathrm{D}^{-} = -\dd_{\bar{\theta}} +\frac{1}{2}\theta
\dd\ ,\ \mathrm{D}^{+} = -\dd_{\theta}
+\frac{1}{2}\bar{\theta}\dd
$$
The $\mathrm{L}$-operator is $s\ell(2|1)$-invariant by
construction and as consequence one obtains the equality
\be \mathbf{M}^{-1}\cdot\mathrm{L}(u)\cdot\mathbf{M} =
\S^{-1}\cdot\mathrm{L}(u)\cdot\S\ ;\ \S
=\mathrm{e}^{\alpha\mathrm{V}_{-}}\cdot
\mathrm{e}^{\bar{\alpha}\mathrm{W}_{-}}\cdot
\mathrm{e}^{\left(\lambda+
\frac{\alpha\bar{\alpha}}{2}\right)\mathrm{S}_{-}}\ ;\
\mathbf{M} = \left(\begin{array}{ccc}
1 & 0 & 0\\
-\bar{\alpha} & 1 & 0\\
\lambda+\frac{\alpha\bar{\alpha}}{2} & -\alpha & 1
   \end{array}\right)
\label{sl21} \ee Finally we put $\Lambda_1 = \left(-\half,
-\half\right)$ in Yang-Baxter equation
$$
\R_{\Lambda_1\Lambda_2}(u-v)\R_{\Lambda_1\Lambda_3}(u)
\R_{\Lambda_2\Lambda_3}(v)=
\R_{\Lambda_2\Lambda_3}(v)\R_{\Lambda_1\Lambda_3}(u)
\R_{\Lambda_1\Lambda_2}(u-v)
$$
change the numeration of the representation spaces
$\Lambda_2 \to \Lambda_1 = (\ell_1, b_1)\ ;\ \Lambda_3 \to
\Lambda_2 = (\ell_2, b_2)$ and consider the restriction on
the invariant subspace
$\mathrm{V}\otimes\mathrm{V}_{\Lambda_1}
\otimes\mathrm{V}_{\Lambda_2}$. In this way one obtains the
defining equation for the $\R$-operator
$$
\mathrm{L}_1(u-v)\mathrm{L}_2(u)
\R_{\Lambda_1\Lambda_2}(v)=
\R_{\Lambda_1\Lambda_2}(v)\mathrm{L}_2(u)\mathrm{L}_1(u-v)
$$
The operator $\mathrm{L}_k$ acts nontrivially on the tensor
product $\mathrm{V}\otimes\mathrm{V}_{\Lambda_k}$ which is
isomorphic to $\mathrm{V}\otimes\C[Z_k]$ where $Z_k
=(z_k,\theta_k,\bar{\theta}_k)$ and the operator
$\R_{\Lambda_1\Lambda_2}(u)$ acts nontrivially on the
tensor product
$\mathrm{V}_{\Lambda_1}\otimes\mathrm{V}_{\Lambda_2}$ which
is isomorphic to $\C[Z_1,Z_2] = \C[Z_1]\otimes \C[Z_2]$.
Note that obtained defining equation is slightly different
from the ones which was used in~\cite{I} and~\cite{DKK}.
The defining equation which is similar to~\cite{I,DKK} is
\be \R^{-1}_{\Lambda_1\Lambda_2}(v-u)
\mathrm{L}_1(u)\mathrm{L}_2(v)=
\mathrm{L}_2(v)\mathrm{L}_1(u)
\R^{-1}_{\Lambda_1\Lambda_2}(v-u) \label{def21} \ee There
exists the well known automorphism of the Yang-Baxter
equation $\R_{\Lambda_1\Lambda_2}(u) \to
\R^{-1}_{\Lambda_1\Lambda_2}(-u)$. In the simplest
$s\ell(2)$ case we have $\R_{\ell_1\ell_2}(u) \sim
\R^{-1}_{\ell_1\ell_2}(-u)$ but for the more complicated
algebras the action of this automorphism is nontrivial. To
proceed in close analogy with~\cite{I,DKK} we shall use the
defining equation~(\ref{def21}) so that we derive the
expression for the operator
$\R^{-1}_{\Lambda_1\Lambda_2}(v-u)$.

\section{The general R-matrix}
\setcounter{equation}{0}

It is useful to extract the operator of permutation
$$
\P_{12}: \C[Z_1]\otimes\C[Z_2] \to \C[Z_2]\otimes\C[Z_1]\
;\ \P_{12}\Psi\left(Z_1 , Z_2\right) = \Psi\left(Z_2 ,
Z_1\right)
$$
from the $\R$-operator $\R^{-1}_{\Lambda_1\Lambda_2}(v-u) =
\P_{12}\check{\R}_{\Lambda_1\Lambda_2}(u;v)$ and solve the
defining equation for the $\check{\R}$-operator. The main
defining equation for the $\check{\R}$-operator is
$$
\check{\R}(u;v)
\mathrm{L}_{1}(u_1,u_2,u_3)\mathrm{L}_{2}(v_1,v_2,v_3)=
\mathrm{L}_{1}(v_1,v_2,v_3)\mathrm{L}_{2}(u_1,u_2,u_3)
\check{\R}(u;v)
$$
$$
u_1 = u+b_1+\ell_1\ ,\ u_2 = u+2b_1\ ,\ u_3 = u+b_1-\ell_1
\ ;\ v_1 = v+b_2+\ell_2\ ,\ v_2 = v+2b_2\ ,\ v_3 =
v+b_2-\ell_2
$$
The operator $\check{\R}$ interchanges all parameters in
the product of two $\mathrm{L}$-operators and similar to
the $s\ell(3)$-case $\check{\R}$-operator can be
represented as the product of the simpler "elementary
building blocks" - $\F$-operators.
\begin{prop}
There exists operator $\F_{1}$ which is the solution of the
defining equations
\begin{equation}
\F_{1}
\mathrm{L}_{1}(u_1,u_2,u_3)\mathrm{L}_{2}(v_1,v_2,v_3)=
\mathrm{L}_{1}(v_1,u_2,u_3)\mathrm{L}_{2}(u_1,v_2,v_3)\F_{1}
\label{sF1}
\end{equation}
$$
\F_1 = \F_1(u_1|v_1,v_2,v_3)\ ;\ \F_1(u_1|v_1,v_2,v_3) =
\F_1(u_1+\lambda|v_1+\lambda,v_2+\lambda,v_3+\lambda)
$$
and these requirements fix the operator $\F_1$ up to
overall normalization constant
$$
\F_1 \sim \S_1^{-1}\cdot\left[
\frac{\Gamma(z_2\partial_2+u_1-v_3+1)}
{\Gamma(z_2\partial_2+v_1-v_3+1)}
\left(f_1+\bar\theta_2\partial_{\bar\theta_2}
\right)-\frac{\Gamma(z_2\partial_2+u_1-v_3)}
{\Gamma(z_2\partial_2+v_1-v_3+1)}
z_2\partial_{\theta_2}\partial_{\bar\theta_2}\right]
\cdot\S_1
$$
$$
\S_1 =
\mathrm{e}^{\frac{\theta_2\bar\theta_2}{2}\partial_2}\cdot
\mathrm{e}^{\theta_1 V_2^{-}}\cdot\mathrm{e}^{\bar\theta_1
W_2^{-}}\cdot
\mathrm{e}^{-\left(z_1+\frac{\theta_1\bar\theta_1}{2}\right)
S_2^{-}}\ ;\ f_1 = \frac{v_1-v_2}{u_1-v_1}
$$
\end{prop}

\begin{prop}
There exists operator $\F_{2}$ which is the solution of the
defining equations
\begin{equation}
\F_{2}
\mathrm{L}_{1}(u_1,u_2,u_3)\mathrm{L}_{2}(v_1,v_2,v_3)=
\mathrm{L}_{1}(u_1,v_2,u_3)\mathrm{L}_{2}(v_1,u_2,v_3)\F_{2}
\label{sF2}
\end{equation}
$$
\F_2 = \F_2(u_1,u_2|v_2,v_3)\ ;\ \F_2(u_1,u_2|v_2,v_3) =
\F_1(u_1+\lambda,u_2+\lambda|v_2+\lambda,v_3+\lambda)
$$
and these requirements fix the operator $\F_2$ up to
overall normalization constant
$$
\F_2 \sim \S_2^{-1}\cdot\left[ f_2 +
u_{12}\cdot\theta_{2}\partial_{\theta_2}+
v_{23}\cdot\bar\theta_{1}\partial_{\bar\theta_1}+
(z_{12}+\theta_1\bar\theta_2)
\partial_{\bar\theta_1}\partial_{\theta_2}
+(u_2-v_2)\theta_{2}\bar\theta_{1}
\partial_{\bar\theta_1}\partial_{\theta_2}\right]\cdot
\S_2
$$
$$
\S_2 = \mathrm{e}^{\theta_1\partial_{\theta_2}}\cdot
\mathrm{e}^{\bar\theta_2\partial_{\bar\theta_1}}
\cdot\mathrm{e}^{\frac{\theta_1\bar\theta_1}{2}\partial_1}\cdot
\mathrm{e}^{-\frac{\theta_2\bar\theta_2}{2}\partial_2} \ ;\
f_2 = \frac{u_{21}v_{23}}{v_2-u_2}\ ,\ u_{12} = u_1-u_2\ ,\
v_{23} = v_2-v_3
$$

\end{prop}

\begin{prop}
There exists operator $\F_{3}$ which is the solution of the
defining equations
\begin{equation}
\F_{3}
\mathrm{L}_{1}(u_1,u_2,u_3)\mathrm{L}_{2}(v_1,v_2,v_3)=
\mathrm{L}_{1}(u_1,u_2,v_3)\mathrm{L}_{2}(v_1,v_2,u_3)
\F_{3} \label{sF3}
\end{equation}
$$
\F_3 = \F_3(u_1,u_2,u_3|v_3)\ ;\ \F_3(u_1,u_2,u_3|v_3) =
\F_3(u_1+\lambda,u_2+\lambda,u_3+\lambda|v_3+\lambda)
$$
and these requirements fix the operator $\F_3$ up to
overall normalization constant
$$
\F_3 \sim \S_3^{-1}\cdot\left[
\frac{\Gamma(z_1\partial_1+u_1-v_3+1)}
{\Gamma(z_1\partial_1+u_1-u_3+1)}
\left(f_3+\theta_1\partial_{\theta_1}
\right)+\frac{\Gamma(z_1\partial_1+u_1-v_3)}
{\Gamma(z_1\partial_1+u_1-u_3+1)}
z_1\partial_{\theta_1}\partial_{\bar\theta_1}\right]
\cdot\S_3
$$
$$
\S_3 =
\mathrm{e}^{-\frac{\theta_1\bar\theta_1}{2}\partial_1}\cdot
\mathrm{e}^{\theta_2 V_1^{-}}\cdot\mathrm{e}^{\bar\theta_2
W_1^{-}}\cdot
\mathrm{e}^{-\left(z_2+\frac{\theta_2\bar\theta_2}{2}\right)
S_1^{-}}\ ;\ f_3 = \frac{u_2-u_3}{u_3-v_3}
$$
\end{prop}

\begin{prop}
The $\check{\R}$-operator can be factorized as follows
$$
\check{\R}(u;v)= \F_1(u_1;v_1,u_2,u_3)
\F_2(u_1,u_2;v_2,u_3)\F_3(u_1,u_2,u_3;v_3)
$$
\end{prop}
There exist six equivalent ways to represent $\check{\R}$
in an factorized form which differ by the order of
$\F$-operators and their parameters. All these expressions
and the proof of the factorization of the
$\check{\R}$-operator can be obtained using the pictures
similar to~\cite{I}.

The defining system of equations for the $\F$-operator can
be reduced to the simpler system which clearly shows the
property of $s\ell(2|1)$-covariance of the $\F$-operator.

\begin{lem}
The defining equation~(\ref{sF1}) for the operator $\F_1$
is equivalent to the system of equations
\begin{equation}
\F_{1} \left[\mathrm{L}_{1}(u_1,u_2,u_3)+
\mathrm{L}_{2}(v_1,v_2,v_3)\right]=
\left[\mathrm{L}_{1}(v_1,u_2,u_3)+
\mathrm{L}_{2}(u_1,v_2,v_3)\right] \F_{1} \label{sF1def}
\end{equation}
$$
\F_{1} z_1 = z_1 \F_{1}\ ,\ \F_{1} \theta_1 = \theta_1
\F_{1}\ ,\ \F_{1} \bar{\theta}_1 = \bar{\theta}_1 \F_{1}
$$
\be\F_1\cdot\left( \mathrm{V}_2^{-}+\bar\theta_1
\mathrm{S}_2^{-}\right) = \left(
\mathrm{V}_2^{-}+\bar\theta_1
\mathrm{S}_2^{-}\right)\cdot\F_1 \label{F1last}
\end{equation}
\end{lem}

\begin{lem}
The defining equation~(\ref{sF2}) for the operator $\F_2$
is equivalent to the system of equations
\begin{equation}
\F_{2} \left[\mathrm{L}_{1}(u_1,u_2,u_3)+
\mathrm{L}_{2}(v_1,v_2,v_3)\right]=
\left[\mathrm{L}_{1}(u_1,v_2,u_3)+
\mathrm{L}_{2}(v_1,u_2,v_3)\right] \F_{2} \label{sF2def}
\end{equation}
$$
\left[\F_{2},z_1-\frac{\theta_1\bar\theta_1}{2}\right] = 0\
,\ \F_{2}\theta_1=\theta_1\F_{2} \ ;
\left[\F_{2},z_2+\frac{\theta_2\bar\theta_2}{2}\right] = 0\
,\ \F_{2} \bar\theta_2 = \bar\theta_2 \F_{2}
$$
\end{lem}

\begin{lem}
The defining equation~(\ref{sF3}) for the operator $\F_3$
is equivalent to the system of equations
\begin{equation}
\F_{3} \left[\mathrm{L}_{1}(u_1,u_2,u_3)+
\mathrm{L}_{2}(v_1,v_2,v_3)\right]=
\left[\mathrm{L}_{1}(u_1,u_2,v_3)+
\mathrm{L}_{2}(v_1,v_2,u_3)\right] \F_{3} \label{sF3def}
\end{equation}
$$
\F_{3} z_2 = z_2 \F_{3}\ ,\ \F_{3} \theta_2 = \theta_2
\F_{3}\ ,\ \F_{3} \bar{\theta}_2 = \bar{\theta}_2 \F_{3}
$$
\be\F_3\cdot\left( \mathrm{W}_1^{-}+\theta_2
\mathrm{S}_1^{-}\right) = \left( \mathrm{W}_1^{-}+\theta_2
\mathrm{S}_1^{-}\right)\cdot\F_3 \label{F3last}
\end{equation}
\end{lem}
The relations in the first line are simply the rules of
commutation of $\F$-operators with $s\ell(2|1)$-generators
written in a compact form. In explicit notations we have
for $\Lambda_1 = (\ell_1, b_1)$ and $\Lambda_2 = (\ell_2,
b_2)$
$$
\F : V_{\Lambda_1}\otimes V_{\Lambda_2} \to
V_{\Lambda^{\prime}_1}\otimes V_{\Lambda^{\prime}_2}
$$
$$
\F_1\ :\ \Lambda^{\prime}_1 = (\ell_1-\xi_1 , b_1+\xi_1)\
;\ \Lambda^{\prime}_2 = (\ell_2+\xi_1 , b_2-\xi_1) \ ;\
\xi_1 = \frac{u_1-v_1}{2}
$$
$$
\F_2\ :\ \Lambda^{\prime}_1 = (\ell_1 , b_1-\xi_2)\ ;\
\Lambda^{\prime}_2 = (\ell_2 , b_2+\xi_2) \ ;\ \xi_2 =
u_2-v_2
$$
$$
\F_3\ :\ \Lambda^{\prime}_1 = (\ell_1+\xi_3 , b_1+\xi_3)\
;\ \Lambda^{\prime}_2 = (\ell_2-\xi_3 , b_2-\xi_3) \ ;\
\xi_3 = \frac{u_3 - v_3}{2}
$$
The $s\ell(2|1)$-invariance of $\R$-matrix follows directly
from the properties of $\F$-operators so that the general
R-matrix $\R^{-1}_{\Lambda_1 \Lambda_2}(v-u) =
\P_{12}\check{\R}_{\Lambda_1\Lambda_2}(u;v)$ is
automatically $s\ell(2|1)$-invariant. \\
{\bf Proof}  Now we are going to the proof of equivalence
of defining equation~(\ref{sF3}) to the
system~(\ref{sF3def}) and derivation of explicit formula
for the operator $\F_3$ . First we show that the
system~(\ref{sF3def}) is the direct consequence of the
eq.~(\ref{sF3}). Let us make the shift $u_k\to u_k+\lambda\
,\ v_1\to v_1+\mu\ ,\ v_2\to v_2+\nu\ ,\ v_3\to
v_3+\lambda$ in the defining equation ~(\ref{sF3}).The
$\F$-operator is invariant under this shift and
$\mathrm{L}$-operators transform as follows
$$
\mathrm{L}_1\to \mathrm{L}_1+\lambda\cdot\II \ ;\
\mathrm{L}_2\to \mathrm{L}_2 +
\lambda\cdot\II+(\mu-\lambda)\left(\begin{array}{ccc}
1 & 0 & 0\\
-\bar\theta_2 & 0 & 0\\
z_2 +\frac{\theta_2\bar\theta_2}{2} & 0 & 0
\end{array}\right) + (\nu-\lambda)\left(\begin{array}{ccc}
0 & 0 & 0\\
\bar\theta_2 & 1 & 0\\
-\theta_2\bar\theta_2 & -\theta_2 & 0
\end{array}\right)
$$
After all one obtains the equation which contains the
arbitrary parameters $\lambda$ , $\mu$ and $\nu$ and as
consequence we derive the system~(\ref{sF3def}) and
equation~(\ref{F3last}). Next we show that from the systems
of equations~(\ref{sF3def}),~(\ref{F3last}) follows
eq.~(\ref{sF3}). This will be almost evident if we rewrite
these equations in equivalent form using the
$s\ell(2|1)$-invariance of the $\mathrm{L}$-operator and
the commutativity of $\F_3$ and
$z_2,\theta_2,\bar\theta_2$. We substitute the factorized
representation~(\ref{factor21}) for the operator
$\mathrm{L}_2$ in the defining equation for the operator
$\F_{3}$
$$
\F_{3} \mathrm{L}_1(u_1,u_2,u_3)\mathbf{M}
   \left(\begin{array}{ccc}
v_1 & \mathrm{D}_2^{-} & -\dd_2 \\
0 & v_2-1 & -\mathrm{D}_2^{+}\\
0 & 0 & v_3
   \end{array}\right)
\mathbf{M}^{-1} = \mathrm{L}_1(u_1,u_2,v_3)\mathbf{M}
   \left(\begin{array}{ccc}
v_1 & \mathrm{D}_2^{-} & -\dd_2 \\
0 & v_2-1 & -\mathrm{D}_2^{+}\\
0 & 0 & u_3
   \end{array}\right)
\mathbf{M}^{-1}\F_{3}
$$
and perform the similarity transformation
$\mathbf{M}^{-1}\cdots\mathbf{M}$ of this matrix equation
using the commutativity $\F_{3}$ and $z_2,
\theta_2,\bar{\theta}_2$. Then using the
$s\ell(2|1)$-invariance of
$\mathrm{L}$-operator~(\ref{sl21})
$$
\mathbf{M}^{-1}\cdot\mathrm{L}_1\cdot\mathbf{M} =
\S^{-1}\cdot \mathrm{L}_1\cdot \S\ ;\ \S =
\mathrm{e}^{\theta_2
\mathrm{V}_1^{-}}\cdot\mathrm{e}^{\bar{\theta}_2
\mathrm{W}_1^{-}}\cdot\mathrm{e}^{-\left(
z_2+\frac{\theta_2\bar{\theta}_2}{2}\right)\mathrm{S}_1^{-}}
\ ;\ \mathbf{M} = \left(\begin{array}{ccc}
1 & 0 & 0\\
-\bar{\theta}_2 & 1 & 0\\
z_2+\frac{\theta_2\bar{\theta}_2}{2} & -\theta_2 & 1
   \end{array}\right)
$$
we derive the equation for the transformed operator
$\mathbf{R} = \S\cdot\F_3\cdot\S^{-1}$ \be \mathbf{R}\cdot
\mathrm{L}_1(u_1,u_2,u_3)\mathbf{L}(v_1,v_2,v_3) =
\mathrm{L}_1(u_1,u_2,v_3)\mathbf{L}(v_1,v_2,u_3)
\cdot\mathbf{R} \label{defF} \ee where
$$
\mathbf{L}(v_1,v_2,v_3)\equiv \S
\cdot\left(\begin{array}{ccc}
v_1 & \mathrm{D}_2^{-} & -\dd_2 \\
0 & v_2-1 & -\mathrm{D}_2^{+}\\
0 & 0 & v_3
   \end{array}\right)\cdot\S^{-1} =
\left(\begin{array}{ccc}
v_1 & \mathrm{D}_2^{-}+\mathrm{W}_{1}^{-}& -\dd_2 +\dd_1\\
0 & v_2-1 & -\mathrm{D}_2^{+}-\mathrm{V}_1^{-}\\
0 & 0 & v_3\end{array}\right)
$$
To derive the system of equations which is equivalent to
the system~(\ref{sF3def}),~(\ref{F3last}) written in terms
of $\mathbf{R}$ we repeat the same trick with the shift of
parameters and obtain the system of equations \be
\mathbf{R}\cdot
\left[\mathrm{L}_1(u_1,u_2,u_3)+\mathbf{L}(v_1,v_2,v_3)\right]
= \left[\mathrm{L}_1(u_1,u_2,v_3)+\mathbf{L}(v_1,v_2,u_3)
\right]\cdot\mathbf{R}\label{1+1} \ee
\begin{equation}
\mathbf{R}\cdot \mathrm{L}_1(u_1,u_2,u_3)
   \left(\begin{array}{ccc}
1 & 0 \\
0 & 1 \\
0 & 0
   \end{array}\right)
= \mathrm{L}_1(u_1,u_2,v_3)
   \left(\begin{array}{ccc}
1 & 0 \\
0 & 1 \\
0 & 0
   \end{array}\right)\mathbf{R}
\label{11}
\end{equation}
It is evident that all equations of the system~(\ref{11})
contained in the equation~(\ref{1+1}) except only one
$(12)$-equation $\mathbf{R} W_1^{-} = W_1^{-} \mathbf{R}$.
We use the system of equation \be \mathbf{R}\cdot
\left[\mathrm{L}_1(u_1,u_2,u_3)+\mathbf{L}(v_1,v_2,v_3)\right]
= \left[\mathrm{L}_1(u_1,u_2,v_3)+\mathbf{L}(v_1,v_2,u_3)
\right]\cdot\mathbf{R}\ ;\  \mathbf{R} W_1^{-} = W_1^{-}
\mathbf{R} \label{DefF} \ee as defining system for operator
$\mathbf{R}$. This system is equivalent to the
system~(\ref{sF3def}),~(\ref{F3last}). Returning to the
system~(\ref{defF}) (it is the system~(\ref{sF3}) written
in terms of $\mathbf{R}$) we note that it is possible to
factorize the matrix $diag(v_1\ ; v_2-1\ ;\ 1)$ from the
right
$$
\mathbf{R}\cdot
\mathrm{L}_1(u_1,u_2,u_3)\left(\begin{array}{ccc} 1 &
\frac{\mathrm{D}_2^{-}+\mathrm{W}_{1}^{-}}{v_2-1} &
-\dd_2 +\dd_1\\
0 & 1 & -\mathrm{D}_2^{+}-\mathrm{V}_1^{-}\\
0 & 0 & v_3\end{array}\right) =
\mathrm{L}_1(u_1,u_2,v_3)\left(\begin{array}{ccc} 1 &
\frac{\mathrm{D}_2^{-}+\mathrm{W}_{1}^{-}}{v_2-1} &
-\dd_2 +\dd_1\\
0 & 1 & -\mathrm{D}_2^{+}-\mathrm{V}_1^{-}\\
0 & 0 & u_3\end{array}\right)\cdot\mathbf{R}
$$
In comparison with~(\ref{DefF}) there are three new
equations only
\begin{equation}
\mathbf{R}\cdot \mathrm{L}_1(u_1,u_2,u_3)
   \left(\begin{array}{ccc}
\dd_1\\
-V_1^{-}\\
v_3
   \end{array}\right)
= \mathrm{L}_1(u_1,u_2,v_3)
   \left(\begin{array}{ccc}
\dd_1\\
-V_1^{-}\\
u_3
   \end{array}\right)\mathbf{R}
\end{equation}
Indeed the system~(\ref{DefF}) contains the equations
$[\mathbf{R} , \mathrm{D}_2^{\pm}]=[\mathbf{R} ,\partial_2]
=[\mathbf{R}, \mathrm{W}_1^{-}]= 0$ and by
conditions~(\ref{11}) we obtain the three new equations. It
is easy to check that these equations follow from the
system~(\ref{DefF}). Finally the systems of
equations~(\ref{DefF}) is defining and it remains to find
the solution. First of all $[\mathbf{R} , z_2]=[\mathbf{R}
, \theta_2]=[\mathbf{R} , \bar\theta_2]=[\mathbf{R} ,
\mathrm{D}_2^{\pm}]=[\mathbf{R} ,\partial_2]=0$ and
therefore the operator $\mathbf{R}$ depends on the
variables $z_1,\theta_1,\bar\theta_1$ only. For simplicity
we use the natural transformation
$$
\mathbf{R} =
\mathrm{e}^{\frac{1}{2}\theta_1\bar{\theta}_1\dd_1}\mathbf{r}
\mathrm{e}^{-\frac{1}{2}\theta_1\bar{\theta}_1\dd_1}
$$
change $z_1,\theta_1, \bar{\theta}_1 \to
z,\theta,\bar{\theta}$ and obtain the system of equations
\begin{equation} \mathbf{r}\dd_{\bar{\theta}} =
\dd_{\bar{\theta}}\mathbf{r}\ ;\ \mathbf{r}\left(z\dd
+\bar{\theta}\dd_{\bar{\theta}}\right) = \left(z\dd
+\bar{\theta}\dd_{\bar{\theta}}\right)\mathbf{r}\ ;\
\mathbf{r}\left(\theta\dd_{\theta}-
\bar{\theta}\dd_{\bar{\theta}}\right) =
\left(\theta\dd_{\theta}-
\bar{\theta}\dd_{\bar{\theta}}\right) \mathbf{r} \label{s1}
\end{equation}
$$
\mathbf{r}\left(z^2\dd +
z\left(\theta\dd_{\theta}+\bar{\theta}\dd_{\bar{\theta}}\right)+
z(u_1-u_3)+\theta\bar{\theta}(u_3-u_2)\right) =
$$
\be =\left(z^2\dd +
z\left(\theta\dd_{\theta}+\bar{\theta}\dd_{\bar{\theta}}\right)+
z(u_1-v_3)+\theta\bar{\theta}(v_3-u_2)\right) \mathbf{r}
\label{s2}
\end{equation}
\begin{equation}
\mathbf{r}\left(z\dd_{\bar{\theta}}+ (u_2-u_3)\theta\right)
= \left(z\dd_{\bar{\theta}}+ (u_2-v_3)\theta\right)
\mathbf{r} \label{s3}
\end{equation}
\begin{equation}
\mathbf{r}\left(-z\left(\dd_{\theta}+\bar{\theta}\dd\right)
+\theta\bar{\theta}\dd_{\theta}+
(u_2-u_1)\bar{\theta}\right) =
\left(-z\left(\dd_{\theta}+\bar{\theta}\dd\right)
+\theta\bar{\theta}\dd_{\theta}+
(u_2-u_1)\bar{\theta}\right) \mathbf{r} \label{s4}
\end{equation}
First of all the equation~(\ref{s4}) is not independent. It
is the consequence of equations~(\ref{s2}) and
$\mathbf{r}\dd_{\bar{\theta}} =
\dd_{\bar{\theta}}\mathbf{r}$ due to commutation relation
$[\mathbf{S}^{+},\mathbf{W}^{-}] = \mathbf{W}^{+}$. The
general solution of the equations~(\ref{s1}) is
$$
\mathbf{r} = \mathbf{a}[z\dd]+ \mathbf{b}[z\dd]\cdot
\theta\dd_{\theta}+\mathbf{c}[z\dd]\cdot
z\dd_{\theta}\dd_{\bar{\theta}}
$$
The equations~(\ref{s2}) and ~(\ref{s3}) results in
recurrence relations
$$
\mathbf{a}[z\dd] - \mathbf{a}[z\dd-1] =
(u_2-u_3)\cdot\mathbf{c}[z\dd]\ ;\ \mathbf{b}[z\dd] =
\frac{u_3-v_3}{u_2-u_3}\cdot\mathbf{a}[z\dd]
$$
$$
\mathbf{c}[z\dd+1]\cdot(z\dd+u_1-u_3+1) =
(z\dd+u_1-v_3)\cdot\mathbf{c}[z\dd]
$$
$$
\mathbf{a}[z\dd+1]\cdot(z\dd+u_1-u_3)+
(u_2-u_3)\cdot\mathbf{c}[z\dd+1] = (z\dd+u_1-v_3)\cdot
\mathbf{a}[z\dd]
$$
$$
\mathbf{a}[z\dd] +
\mathbf{b}[z\dd]\cdot(z\dd+u_1-u_3)-(u_2-u_3)\cdot\mathbf{c}[z\dd]
= \mathbf{a}[z\dd-1] +
(z\dd+u_1-v_3)\cdot\mathbf{b}[z\dd-1]
$$
The solution of these equations has the form
$$\mathbf{a}[z\dd] =
\frac{\Gamma(z\dd+u_1-v_3+1)}{\Gamma(z\dd+u_1-u_3+1)}\ ;\
\mathbf{b}[z\dd] = \frac{u_3-v_3}{u_2-u_3}
\frac{\Gamma(z\dd+u_1-v_3+1)}{\Gamma(z\dd+u_1-u_3+1)}
$$
$$
\mathbf{c}[z\dd] = \frac{u_3-v_3}{u_2-u_3}
\frac{\Gamma(z\dd+u_1-v_3)}{\Gamma(z\dd+u_1-u_3+1)}
$$
Collect everything together we obtain the expression for
the operator $\F_3$ from the Proposition. All calculations
for the operator $\F_1$ are very similar.

It remains to prove the equivalence of defining
equation~(\ref{sF2}) to the system~(\ref{sF2def}) and
derive the explicit formula for the operator $\F_2$. First
we show that the system~(\ref{sF2def}) is the direct
consequence of the eq.~(\ref{sF2}). Let us make the shift
$u_1\to u_1+\lambda\ ,\ u_2\to u_2+\lambda\ ,\ u_3\to
u_3+\mu\ ,\ v_1\to v_1+\nu\ ,\ v_2\to v_2+\lambda\ ,\
 v_3\to v_3+\lambda$ in the defining equation ~(\ref{sF2})
for the operator $\F_2$. The $\F$-operator is invariant
under this shift and $\mathrm{L}$-operators transform as
follows
$$
\mathrm{L}_1\to \mathrm{L}_1+ \lambda\cdot\II
+(\mu-\lambda)\left(\begin{array}{ccc}
0 & 0 & 0\\
0 & 0 & 0\\
-z_1 +\frac{\theta_1\bar\theta_1}{2} & \theta_1 & 1
\end{array}\right)\ ;\
\mathrm{L}_2 \to \mathrm{L}_2 + \lambda\cdot\II+
(\nu-\lambda)\left(\begin{array}{ccc}
1 & 0 & 0\\
-\bar\theta_2 & \lambda & 0\\
z_2 +\frac{\theta_2\bar\theta_2}{2} & 0
&0\end{array}\right)
$$
After all one obtains the equation which contains the
arbitrary parameters $\lambda$ , $\mu$ and $\nu$ and as
consequence we derive the system~(\ref{sF2def}). Next we
show that from the systems of equations~(\ref{sF2def})
follows eq.~(\ref{sF2}). This will be almost evident if we
rewrite these equations in equivalent form using the
$s\ell(2|1)$-invariance of the $\mathrm{L}$-operator and
the commutativity of $\F_2$ and
$z_1-\frac{\theta_1\bar\theta_1}{2},\theta_1
,z_2+\frac{\theta_2\bar\theta_2}{2},\bar\theta_2$. First of
all it is useful to make the transformation
$$
\F_2 = \S^{-1}\cdot\mathbf{R}\cdot\S\ ;\ \S
=\mathrm{e}^{\frac{\theta_1\bar\theta_1}{2}\partial_1}\cdot
\mathrm{e}^{-\frac{\theta_2\bar\theta_2}{2}\partial_2}
$$
so that $\mathbf{R}$ commutes with $z_1,\theta_1
,z_2,\bar\theta_2$ now. The corresponding transformation
for the $\mathrm{L}$-operators can be easily derived using
factorized representation~(\ref{factor21}). Next it is
possible to make the two similarity transformations of the
defining equation~(\ref{sF2}) using simple matrices which
commute with operator $\mathbf{F}$. After all these
transformations the defining equation~(\ref{sF2}) for the
$\mathbf{F}$-operator in factorized form looks as follows
$$
\mathbf{R}\cdot\mathbf{l}_1(u_1,u_2,u_3)\cdot
\mathbf{M}\cdot\mathbf{l}_2(v_1,v_2,v_3) =
\mathbf{l}_1(u_1,v_2,u_3)\cdot
\mathbf{M}\cdot\mathbf{l}_2(v_1,u_2,v_3) \cdot\mathbf{R}\
;\ \mathbf{M}\equiv \left(\begin{array}{ccc}
1 & 0 & 0 \\
-\bar\theta_2 & 1 & 0 \\
-z_{12}-\theta_1\bar\theta_2 & \theta_1 & 1
   \end{array}\right)
$$
$$
\mathbf{l}_1(u_1,u_2,u_3)\equiv\left(\begin{array}{ccc}
1 & 0 & 0 \\
-\bar{\theta}_1 & 1 & 0 \\
0 & 0 & 1
   \end{array}\right)\left(\begin{array}{ccc}
u_1 & -\dd_{\bar{\theta}_1} & -\dd_1 \\
0 & u_2-1 & \dd_{\theta_1}-\bar\theta_1\dd_1 \\
0 & 0 & u_3
   \end{array}\right)\left(\begin{array}{ccc}
1 & 0 & 0 \\
\bar{\theta}_1 & 1 & 0 \\
0 & 0 & 1
   \end{array}\right)
$$
$$
\mathbf{l}_2(v_1,v_2,v_3)\equiv \left(\begin{array}{ccc}
1 & 0 & 0 \\
0 & 1 & 0 \\
0 & -\theta_2 & 1
   \end{array}\right)\left(\begin{array}{ccc}
v_1 & -\dd_{\bar{\theta}_2}+\theta_2\dd_2 & -\dd_2 \\
0 & v_2-1 & \dd_{\theta_2} \\
0 & 0 & v_3
   \end{array}\right)\left(\begin{array}{ccc}
1 & 0 & 0 \\
0 & 1 & 0 \\
0 & \theta_2 & 1
   \end{array}\right)
$$
Next step we rewrite the defining equation for the
transformed operator
$$
\mathbf{r} = \mathrm{e}^{\theta_1\partial_{\theta_2}}\cdot
\mathrm{e}^{\bar\theta_2\partial_{\bar\theta_1}}
\cdot\mathbf{R}\cdot
\mathrm{e}^{-\theta_1\partial_{\theta_2}}\cdot
\mathrm{e}^{-\bar\theta_2\partial_{\bar\theta_1}}
$$
in the form \be
\mathbf{r}\cdot\mathbf{L}_1(u_1,u_2,u_3)\mathbf{m}
\mathbf{L}_2(v_1,v_2,v_3) =
\mathbf{L}_1(u_1,v_2,u_3)\mathbf{m}
\mathbf{L}_2(v_1,u_2,v_3) \cdot\mathbf{r} \label{defF2} \ee
where
$$
\mathbf{L}_1(u_1,u_2,u_3) = \left(\begin{array}{ccc}
u_1-1+\bar\theta_1\dd_{\bar\theta_1} & -\dd_{\bar\theta_1}
&
-\dd_1 \\
(u_2-u_1)\bar\theta_1 &
u_2-1+\bar\theta_1\dd_{\bar\theta_1} &
\dd_{\theta_1}-\dd_{\theta_2}-\bar\theta_2\dd_1 \\
0 & 0 & u_3
   \end{array}\right)
$$
$$
\mathbf{L}_2(v_1,v_2,v_3) = \left(\begin{array}{ccc} v_1 &
-\dd_{\bar\theta_2}+\dd_{\bar\theta_1}+\theta_1\dd_2 &
-\dd_2 \\
0 & v_2-\theta_2\dd_{\theta_2} & -\dd_{\theta_2} \\
0 & (v_3-v_2)\theta_2 & v_3-\theta_2\dd_{\theta_2}
   \end{array}\right)
$$
To derive the system of equations which is equivalent to
the system~(\ref{sF2def}) written in terms of $\mathbf{r}$
we repeat the same trick with the shift of parameters and
obtain \be
\mathbf{r}\cdot\left[\mathbf{L}_1(u_1,u_2,u_3)\mathbf{m}+
\mathbf{m}\mathbf{L}_2(v_1,v_2,v_3)\right] =
\left[\mathbf{L}_1(u_1,v_2,u_3)\mathbf{m}+
\mathbf{m}\mathbf{L}_2(v_1,u_2,v_3)\right]\cdot\mathbf{r}
\label{2+2} \ee This system results in a simple equations
\be \mathbf{r}\left(\dd_{\theta_1}-\bar\theta_1\dd_1\right)
= \left(\dd_{\theta_1}-\bar\theta_1\dd_1\right)\mathbf{r}\
;\ \mathbf{r}\left(\dd_{\bar\theta_2}-\theta_1\dd_2\right)
= \left(\dd_{\bar\theta_2}-\theta_1\dd_2\right)\mathbf{r}\
;\ \mathbf{r}(\dd_1+\dd_2) = (\dd_1+\dd_2)\mathbf{r}
\label{eq1}\ee \be \mathbf{r}\left(\theta_2\dd_{\theta_2}+
(z_{12}+\theta_1\bar\theta_2)\dd_2\right) =
\left(\theta_2\dd_{\theta_2}+
(z_{12}+\theta_1\bar\theta_2)\dd_2\right)\mathbf{r}
\label{eq2}\ee \be
\mathbf{r}\left(\bar\theta_1\dd_{\bar\theta_1}+
(z_{12}+\theta_1\bar\theta_2)\dd_1\right) =
\left(\bar\theta_1\dd_{\bar\theta_1}+
(z_{12}+\theta_1\bar\theta_2)\dd_1\right)\mathbf{r}
\label{eq3}\ee \be
\mathbf{r}\left(\bar\theta_1\dd_{\bar\theta_1}-
\theta_2\dd_{\theta_2}\right) =
\left(\bar\theta_1\dd_{\bar\theta_1}-
\theta_2\dd_{\theta_2}\right)\mathbf{r} \label{eq4} \ee \be
\mathbf{r}\left((z_{12}+\theta_1\bar\theta_2)\dd_{\theta_2}+
(u_2-u_1)\bar\theta_1\right) =
\left((z_{12}+\theta_1\bar\theta_2)\dd_{\theta_2}+
(v_2-u_1)\bar\theta_1\right)\mathbf{r} \label{eq5}\ee \be
\mathbf{r}\left((z_{12}+\theta_1\bar\theta_2)\dd_{\bar\theta_1}+
(v_3-v_2)\theta_2\right) =
\left((z_{12}+\theta_1\bar\theta_2)\dd_{\bar\theta_1}+
(v_3-v_2)\theta_2\right)\mathbf{r} \label{eq6} \ee
Returning to the system~(\ref{defF2})(it is the
system~(\ref{sF2}) written in terms of $\mathbf{r}$) we
note that it is possible to factorize the matrix $diag(v_1
; 1 ; 1)$ from the right and the matrix $diag(1 ; 1 ; u_3)$
from the left so that $v_1,u_3$-dependence disappear from
equation. The obtained system of equations is the
consequence of the system~(\ref{eq1})-(\ref{eq6}). The
proof is by direct calculation. Now we are going to the
solution of the defining system of equations. The general
solution of equations~(\ref{eq1})-(\ref{eq4}) has the form
$$
\mathbf{r} = \mathbf{a}
+\mathbf{b}\cdot\bar\theta_1\dd_{\bar\theta_1}+
\mathbf{c}\cdot\theta_2\dd_{\theta_2} +
\mathbf{d}\cdot(z_{12}+\theta_1\bar\theta_2)
\dd_{\bar\theta_1}\dd_{\theta_2} +
\mathbf{e}\cdot\bar\theta_1\theta_2\cdot
\dd_{\bar\theta_1}\dd_{\theta_2}
$$
where
$\mathbf{a},\mathbf{b},\mathbf{c},\mathbf{d},\mathbf{e}$
are some constants. The equations~(\ref{eq5})
and~(\ref{eq6}) fix these constants
$$
\mathbf{a} =
\frac{(u_2-u_1)(v_2-v_3)}{v_2-u_2}\cdot\mathbf{d}\ ;\
\mathbf{b} = (v_2-v_3)\cdot\mathbf{d}\ ;\ \mathbf{c} =
(u_1-u_2)\cdot\mathbf{d}\ ;\ \mathbf{e} =
(u_2-v_2)\cdot\mathbf{d}
$$
Collect everything together we obtain the expression for
the operator $\F_2$ from the Proposition.

\section{Conclusions}

We have shown that the general R-matrix can be represented
as the product of the simple "building blocks" --
$\F$-operators. In the first paper~\cite{I} we have
demonstrated how this factorization arises in the simplest
situations of the symmetry algebra $s\ell(2)$ and
$s\ell(3)$. In the present paper we have showed that the
same factorization take place for the R-matrix with
supersymmetry algebra $s\ell(2|1)$. It seems that this
phenomenon is quite general and all results can be
generalized to the symmetry algebra $s\ell(n)$ and to the
supersymmetry algebra $s\ell(n|m)$.

\section{Acknowledgments}

I would like to thank D.Karakhanyan, R.Kirschner,
G.Korchemsky, P.Kulish and A.Manashov for the stimulating
discussions and critical remarks on the different stages of
this work. This work was supported by the grant 03-01-00837
of the Russian Foundation for Fundamental Research.

\section*{Appendix }
\setcounter{equation}{0}

In this Appendix we calculate the matrix elements of the
$\F$-operators and as consequence obtain the matrix
elements of R-matrix. It is additional check of the main
results and after all we recover the formulae from the
paper~\cite{DKK}.

All lowest weights in the space $V_{\ell_1,b_1}\otimes
V_{\ell_2,b_2}$ are divided on two sets. There are the even
lowest weights
$$
 \Phi^{\pm}_n\equiv
\left(Z_{12}\pm\half\theta_{12}\bar\theta_{12}\right)^{n} \
;\ D_1^{\pm}\Phi^{\pm}_n = 0 \ ,\ S \Phi^{\pm}_n = (n
+\ell_1+\ell_2)\Phi^{\pm}_n \ ,\ B \Phi^{\pm}_n = (b_1+b_2)
\Phi^{\pm}_n
$$
and the odd lowest weights
$$
\Psi^{-}_n\equiv\theta_{12}Z_{12}^{n} \ ;\
\Psi^{+}_n\equiv\bar\theta_{12}Z_{12}^{n} \ ;\ S
\Psi^{\pm}_n = \left(n
+\ell_1+\ell_2+\frac{1}{2}\right)\Psi^{\pm}_n \ ,\ B
\Psi^{\pm}_n =
\left(b_1+b_2\pm\frac{1}{2}\right)\Psi^{\pm}_n
$$
In this section we shall calculate the
action of $\F$-operators on these lowest weights.

\section*{Operator $\F_3$}

We have factorized representation for the operator $\F_3$
$$
\F_3 \sim \S_3^{-1}\cdot \mathbf{r}_3 \cdot\S_3 \ ;\
\mathbf{r}_3 \equiv \frac{\Gamma(z_1\partial_1+u_1-v_3+1)}
{\Gamma(z_1\partial_1+u_1-u_3+1)}
\left(\frac{u_2-u_3}{u_3-v_3}+\theta_1\partial_{\theta_1}
\right)+\frac{\Gamma(z_1\partial_1+u_1-v_3)}
{\Gamma(z_1\partial_1+u_1-u_3+1)}
z_1\partial_{\theta_1}\partial_{\bar\theta_1}.
$$
In the explicit form the action of the operators $\S_3$ and
$\S_3^{-1}$ is
$$\S_3 \Phi\left(z_1,\theta_1
,\bar\theta_1|z_2,\theta_2 ,\bar\theta_2\right) =
\Phi\left(z_1+z_2+\frac{\theta_2\bar\theta_1}{2}
-\frac{\theta_1\bar\theta_2}{2}-\frac{\theta_1\bar\theta_1}{2},
\theta_1+\theta_2, \bar\theta_1+\bar\theta_2 | z_2,\theta_2
,\bar\theta_2\right)$$$$ \S_3^{-1} \Phi\left(z_1,\theta_1
,\bar\theta_1| z_2,\theta_2 ,\bar\theta_2\right) =
\Phi\left(z_1-z_2+\frac{\theta_1\bar\theta_1}{2}
+\frac{\theta_2\bar\theta_2}{2}-\theta_2\bar\theta_1,
\theta_1-\theta_2, \bar\theta_1-\bar\theta_2 | z_2,\theta_2
,\bar\theta_2\right)$$ First we calculate the action of
$\S_3$ $$\S_3 \Phi_n^{+} = z_1^n \ ;\ \S_3 \Phi_n^{-} =
\left(z_1-\theta_1\bar\theta_1\right)^n
$$ then the action
of $\mathbf{r}_3$
$$ \mathbf{r}_3\cdot z_1^n =
\frac{u_2-u_3}{u_3-v_3}
\frac{\Gamma(n+u_1-v_3+1)}{\Gamma(n+u_1-u_3+1)}\cdot z_1^n
$$
$$
\mathbf{r}_3\cdot \left(z_1-\theta_1\bar\theta_1\right)^n =
\frac{u_2-v_3}{u_3-v_3}
\frac{\Gamma(n+u_1-v_3)}{\Gamma(n+u_1-u_3)}\cdot
\left(z_1-\theta_1\bar\theta_1\right)^n +(u_2-u_1)\cdot
\frac{\Gamma(n+u_1-v_3)}{\Gamma(n+u_1-u_3+1)}\cdot z_1^n
$$
and finally action of $\S_3^{-1}$
$$ \Phi_n^{+} =
\S_3^{-1}\cdot z_1^n \ ;\ \Phi_n^{-} =
\S_3^{-1}\left(z_1-\theta_1\bar\theta_1\right)^n
$$ so that one obtains
$$
\F_3\cdot \Phi_n^{+} = \frac{u_2-u_3}{u_3-v_3}
\frac{\Gamma(n+u_1-v_3+1)}{\Gamma(n+u_1-u_3+1)}\cdot
\Phi_n^{+}
$$
$$
\F_3\cdot \Phi_n^{-} = \frac{u_2-v_3}{u_3-v_3}
\frac{\Gamma(n+u_1-v_3)}{\Gamma(n+u_1-u_3)}\cdot \Phi_n^{-}
+(u_2-u_1)\cdot
\frac{\Gamma(n+u_1-v_3)}{\Gamma(n+u_1-u_3+1)}\cdot
\Phi_n^{+}
$$
For the odd lowest weights all is simpler
$$ \S_3
\Psi_n^{+} = \bar\theta_1\cdot z_1^n \ ;\ \S_3 \Psi_n^{-} =
\theta_1\cdot z_1^n
$$
$$
\mathbf{r}_3\cdot \bar\theta_1 z_1^n =
\frac{u_2-u_3}{u_3-v_3}
\frac{\Gamma(n+u_1-v_3+1)}{\Gamma(n+u_1-u_3+1)}\cdot
\bar\theta_1 z_1^n\ ;\ \mathbf{r}_3\cdot \theta_1 z_1^n =
\frac{u_2-v_3}{u_3-v_3}
\frac{\Gamma(n+u_1-v_3+1)}{\Gamma(n+u_1-u_3+1)}\cdot
\theta_1 z_1^n
$$
and finally we have
$$
\F_3\cdot \Psi_n^{+} = \frac{u_2-u_3}{u_3-v_3}
\frac{\Gamma(n+u_1-v_3+1)}{\Gamma(n+u_1-u_3+1)}\cdot
\Psi_n^{+}\ ;\  \F_3\cdot \Psi_n^{-} =
\frac{u_2-v_3}{u_3-v_3}
\frac{\Gamma(n+u_1-v_3+1)}{\Gamma(n+u_1-u_3+1)}\cdot
\Psi_n^{-}
$$

\section*{Operator $\F_2$}

We have factorized representation for the operator $\F_2$
$$ \F_2 \sim \S^{-1}\cdot \mathbf{r}_2 \cdot\S $$
$$
\mathbf{r}_2 \equiv \frac{(u_2-u_1)(v_2-v_3)}{v_2-u_2}+
(u_2-u_1)\theta_{12}\partial_{\theta_2}+
(v_2-v_3)\bar\theta_{12}\partial_{\bar\theta_1} +
(z_{12}+\theta_1\bar\theta_2)
\partial_{\bar\theta_1}\partial_{\theta_2}
+(v_2-u_2)\theta_{12}\bar\theta_{12}
\partial_{\bar\theta_1}\partial_{\theta_2}
$$
In the explicit form the action of the operators $\S$ and
$\S^{-1}$ is
$$
\S\Phi\left(z_1,\theta_1 ,\bar\theta_1|z_2,\theta_2
,\bar\theta_2\right) =
\Phi\left(z_1+\frac{\theta_1\bar\theta_1}{2},
\theta_1,\bar\theta_1|z_2 -\frac{\theta_2\bar\theta_2}{2},
\theta_2, \bar\theta_2\right)
$$
$$
\S^{-1} \Phi\left(z_1,\theta_1 ,\bar\theta_1|z_2,\theta_2
,\bar\theta_2\right) =
\Phi\left(z_1-\frac{\theta_1\bar\theta_1}{2},
\theta_1,\bar\theta_1|z_2 +\frac{\theta_2\bar\theta_2}{2},
\theta_2, \bar\theta_2\right)
$$
First we calculate the action of $\S$
$$
\S\Phi_n^{+} = \left(z_{12}+\theta_1\bar\theta_2+
\theta_{12}\bar\theta_{12}\right)^n \ ;\ \S\Phi_n^{-} =
\left(z_{12}+\theta_1\bar\theta_2\right)^n
$$ then the
action of $\mathbf{r}_2$
$$
\mathbf{r}_2\cdot \left(z_{12}+\theta_1\bar\theta_2+
\theta_{12}\bar\theta_{12}\right)^n =
\frac{(u_2-v_3)(v_2-u_1)}{v_2-u_2}\cdot
\left(z_{12}+\theta_1\bar\theta_2+
\theta_{12}\bar\theta_{12}\right)^n -(u_1-v_3+n)\cdot
\left(z_{12}+\theta_1\bar\theta_2\right)^n
$$
$$
\mathbf{r}_2\cdot
\left(z_{12}+\theta_1\bar\theta_2\right)^n =
\frac{(u_2-u_1)(v_2-v_3)}{v_2-u_2}\cdot
\left(z_{12}+\theta_1\bar\theta_2\right)^n
$$ and finally one obtains
$$
\F_2\cdot \Phi_n^{+} =
\frac{(u_2-v_3)(v_2-u_1)}{v_2-u_2}\cdot \Phi_n^{+}
-(u_1-v_3+n)\cdot\Phi_n^{-}\ ;\ \F_2\cdot \Phi_n^{-} =
\frac{(u_2-u_1)(v_2-v_3)}{v_2-u_2}\cdot \Phi_n^{-} $$ For
the odd lowest weights all is simpler
$$
\S\Psi_n^{+} = \bar\theta_{12}\cdot
\left(z_{12}+\theta_1\bar\theta_2\right)^n \ ;\
\S\Psi_n^{-} = \theta_{12}\cdot
\left(z_{12}+\theta_1\bar\theta_2\right)^n
$$
$$
\mathbf{r}_2\cdot \bar\theta_{12}\cdot
\left(z_{12}+\theta_1\bar\theta_2\right)^n =
\frac{(v_2-u_1)(v_2-v_3)}{v_2-u_2}\cdot\bar\theta_{12}\cdot
\left(z_{12}+\theta_1\bar\theta_2\right)^n
$$
$$
\mathbf{r}_2\cdot \theta_{12}\cdot
\left(z_{12}+\theta_1\bar\theta_2\right)^n =
\frac{(u_2-u_1)(u_2-v_3)}{v_2-u_2}\cdot\theta_{12}\cdot
\left(z_{12}+\theta_1\bar\theta_2\right)^n
$$ and finally we have
$$
\F_2 \cdot \Psi_n^{+} =
\frac{(v_2-u_1)(v_2-v_3)}{v_2-u_2}\cdot\Psi_n^{+}\ ;\ \F_2
\cdot \Psi_n^{-} =
\frac{(u_2-u_1)(u_2-v_3)}{v_2-u_2}\cdot\Psi_n^{-}.
$$

\section*{Operator $\F_1$}

We have factorized representation for the operator $\F_1$
$$
\F_1 \sim \S_1^{-1}\cdot \mathbf{r}_1 \cdot\S_1 \ ;\
\mathbf{r}_1 \equiv \frac{\Gamma(z_2\partial_2+u_1-v_3+1)}
{\Gamma(z_2\partial_2+v_1-v_3+1)}
\left(\frac{v_1-v_2}{u_1-v_1}+\bar\theta_2\partial_{\bar\theta_2}
\right)-\frac{\Gamma(z_2\partial_2+u_1-v_3)}
{\Gamma(z_2\partial_2+v_1-v_3+1)}
z_2\partial_{\theta_2}\partial_{\bar\theta_2}
$$
and the action of the operators $\S_1$ and $\S_1^{-1}$ in
explicit form is
$$
\S_1 \Phi\left(z_1,\theta_1 ,\bar\theta_1 |z_2,\theta_2
,\bar\theta_2\right) = \Phi\left(z_1,\theta_1 ,\bar\theta_1
|z_2+z_1+\frac{\theta_2\bar\theta_2}{2}
+\frac{\theta_1\bar\theta_2}{2}-\frac{\theta_2\bar\theta_1}{2},
\theta_2+\theta_1, \bar\theta_2+\bar\theta_1 \right)
$$
$$
\S_1^{-1} \Phi\left(z_1,\theta_1 ,\bar\theta_1
|z_2,\theta_2 ,\bar\theta_2\right) = \Phi\left(z_1,\theta_1
,\bar\theta_1 |z_2-z_1-\frac{\theta_1\bar\theta_1}{2}
-\frac{\theta_2\bar\theta_2}{2}+\theta_2\bar\theta_1,
\theta_2-\theta_1, \bar\theta_2-\bar\theta_1 \right)
$$
First we calculate the action of $\S_1$
$$
\S_1 \Phi_n^{+} = (-z_2)^n \ ;\ \S_1 \Phi_n^{-} =
\left(-z_1-\theta_2\bar\theta_2\right)^n
$$ then the action of $\mathbf{r}_1$
$$
\mathbf{r}_1\cdot (-z_2)^n = \frac{v_1-v_2}{u_1-v_1}
\frac{\Gamma(n+u_1-v_3+1)}{\Gamma(n+v_1-v_3+1)}\cdot
(-z_2)^n
$$
$$
\mathbf{r}_1\cdot \left(-z_2-\theta_2\bar\theta_2\right)^n
= \frac{u_1-v_2}{u_1-v_1}
\frac{\Gamma(n+u_1-v_3)}{\Gamma(n+v_1-v_3)}\cdot
\left(-z_2-\theta_2\bar\theta_2\right)^n +(v_3-v_2)\cdot
\frac{\Gamma(n+u_1-v_3)}{\Gamma(n+v_1-v_3+1)}\cdot (-z_2)^n
$$
and finally one obtains
$$
\F_1\cdot \Phi_n^{+} =
\frac{v_1-v_2}{u_1-v_1}
\frac{\Gamma(n+u_1-v_3+1)}{\Gamma(n+v_1-v_3+1)}\cdot
\Phi_n^{+}
$$
$$
\F_1\cdot \Phi_n^{-} = \frac{u_1-v_2}{u_1-v_1}
\frac{\Gamma(n+u_1-v_3)}{\Gamma(n+v_1-v_3)}\cdot \Phi_n^{-}
+(v_3-v_2)\cdot
\frac{\Gamma(n+u_1-v_3)}{\Gamma(n+v_1-v_3+1)}\cdot
\Phi_n^{+}
$$
For the odd lowest weights we have
$$
\F_1\cdot \Psi_n^{+} = \frac{u_1-v_2}{u_1-v_1}
\frac{\Gamma(n+u_1-v_3+1)}{\Gamma(n+v_1-v_3+1)}\cdot
\Psi_n^{+}\ ;\ \F_1\cdot \Psi_n^{-} =
\frac{v_1-v_2}{u_1-v_1}
\frac{\Gamma(n+u_1-v_3+1)}{\Gamma(n+v_1-v_3+1)}\cdot
\Psi_n^{-}
$$

\section*{Operator R}

The matrix elements for the $\R$-operator are obtained from
the matrix elements of $\F$-operators by the formula
$$
\R(u;v)\sim \F_1(u_1;v_1,u_2,u_3)
\F_2(u_1,u_2;v_2,u_3)\F_3(u_1,u_2,u_3;v_3)
$$
The result of calculations is the following
$$
\R\Phi_n^{+} \sim \mathrm{R}\cdot\left\{
(u_2-u_3)(v_2-v_1)\cdot
\frac{\Gamma(n+u_1-v_3+1)}{\Gamma(n+v_1-u_3+1)}\cdot\Phi_n^{+}
+ (u_2-v_2)\cdot
\frac{\Gamma(n+u_1-v_3+1)}{\Gamma(n+v_1-u_3)}\cdot\Phi_n^{+}
\right\}
$$
$$
\R\Phi_n^{-} \sim \mathrm{R}\cdot\left\{
(u_2-u_1)(v_2-v_3)\cdot
\frac{\Gamma(n+u_1-v_3)}{\Gamma(n+v_1-u_3)}\cdot\Phi_n^{-}
+ \mathrm{C}\cdot
\frac{\Gamma(n+u_1-v_3)}{\Gamma(n+v_1-u_3+1)}\cdot\Phi_n^{+}
\right\}
$$
where
$$\mathrm{C} = (u_2-v_3)(v_2-u_3)(u_1-v_1)+
(v_2-u_1)(v_1-u_2)(v_3-u_3)+ (u_1-v_1)(u_2-v_2)(u_3-v_3)
$$
$$
\R \Psi_n^{+} \sim \mathrm{R}\cdot (v_2-u_1)(v_2-u_3)\cdot
\frac{\Gamma(n+u_1-v_3+1)}{\Gamma(n+v_1-u_3+1)}\cdot\Psi_n^{+}
$$
$$
\R \Psi_n^{-} \sim \mathrm{R}\cdot
(u_2-v_1)(u_2-v_3)\cdot
\frac{\Gamma(n+u_1-v_3+1)}{\Gamma(n+v_1-u_3+1)}\cdot\Psi_n^{+}
$$
We extract the common normalization factor
$$
\mathrm{R} \equiv
\frac{(u_2-u_1)(u_2-u_3)}{(u_1-v_1)(u_2-v_2)(u_3-v_3)}.
$$
After substitution of parameters in explicit form
$$
u_1 = u+b_1+\ell_1\ ,\ v_1 = v+b_2+\ell_2 \ ;\ u_2 = u+2
b_1\ ,\ v_2 = v+2 b_2 \ ;\ u_3 = u+b_1-\ell_1 \ ,\ v_3 =
v+b_2-\ell_2
$$
we recover the formulae for the matrix elements of
R-operator from the paper~\cite{DKK}.

\end{document}